\documentclass[11pt]{amsart}
\usepackage{amsmath,amssymb,amsthm,graphicx,geometry,hyperref,placeins}
\hypersetup{colorlinks=true,linkcolor=blue,citecolor=blue,urlcolor=blue}
\allowdisplaybreaks[2]
\theoremstyle{plain}
\newtheorem{thm}{Theorem}
\newtheorem{lem}{Lemma}
\newtheorem{prop}{Proposition}
\theoremstyle{remark}
\newtheorem*{rem}{Remark}
\newcommand{\Pow}{\operatorname{Pow}}

\begin{document}
\title[New Theorems on the Pentagon and the Pentagram]{New Theorems on the Pentagon and the Pentagram}
\author{Quang Hung Tran}
\address{High School for Gifted Students, Vietnam National University, Hanoi, Vietnam}
\email{tranquanghung@hus.edu.vn}
\begin{abstract}
We prove two pencil theorems for circumcenters in Miquel's pentagon and pentagram configuration. Each construction joins five triangle centers to five centers of circles through the Miquel points. An explicit formula expresses the finite common point as an affine combination of the ten triangle centers. Its weights are signed areas divided by squared circumradii and depend only on the side directions. Their sum also detects the parallel case. When either the vertices or the outer vertices are concyclic, the finite pencil point lies on the line joining their circle center to the Miquel circle center. We also iterate the construction when the ear centers and Miquel points lie on one circle. Radical axes give the resulting collinearities, and the common circle configuration is identified with the classical branch of Morley's five circles theorem.
\end{abstract}
\maketitle

\section{Introduction}
\begin{figure}[htbp]
\centering\includegraphics[width=.98\linewidth]{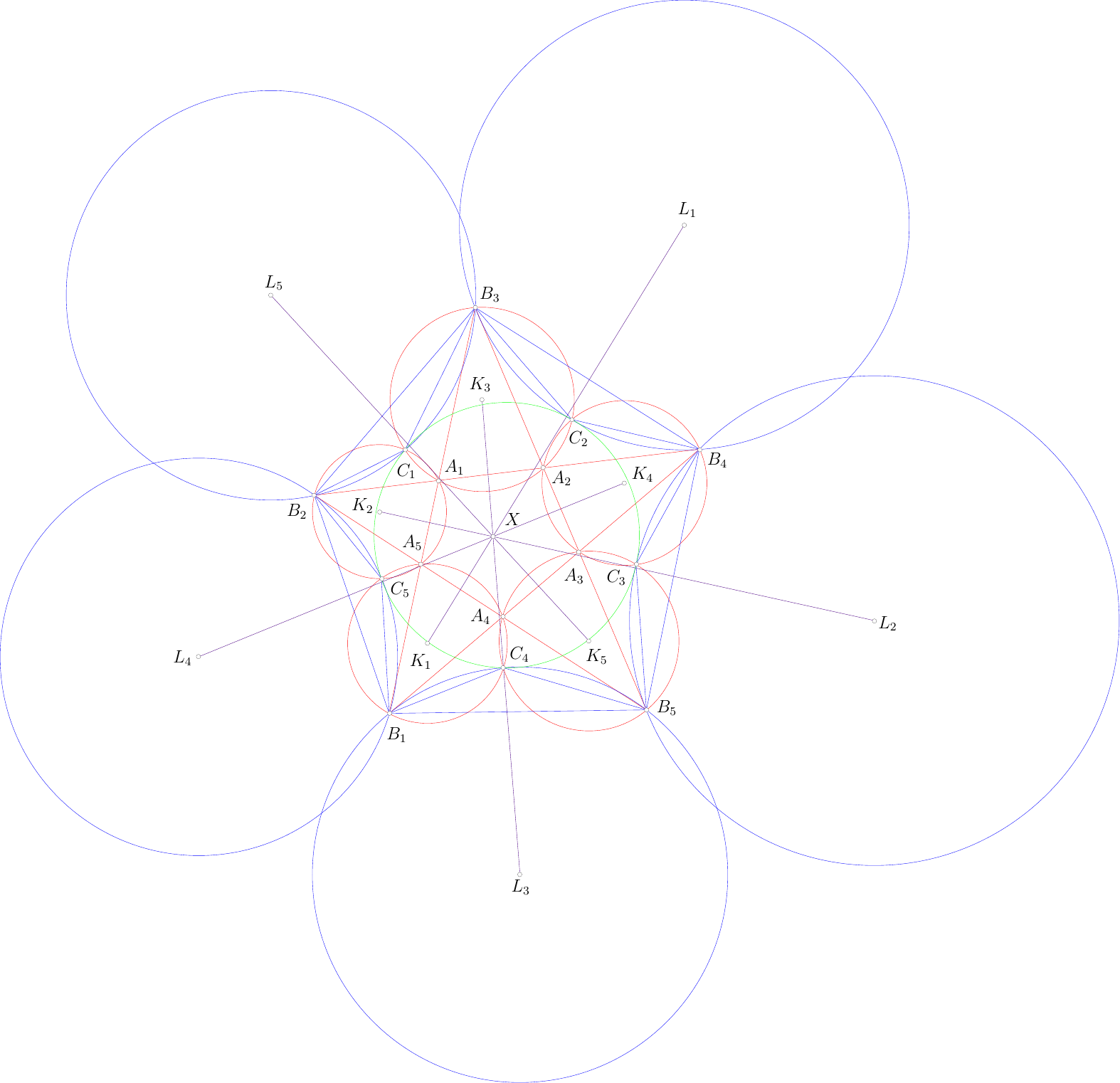}
\caption{The principal circumcenter pencil, here with a finite common point.}\label{fig:pencil}
\end{figure}

\begin{figure}[htbp]
\centering\includegraphics[width=.93\linewidth]{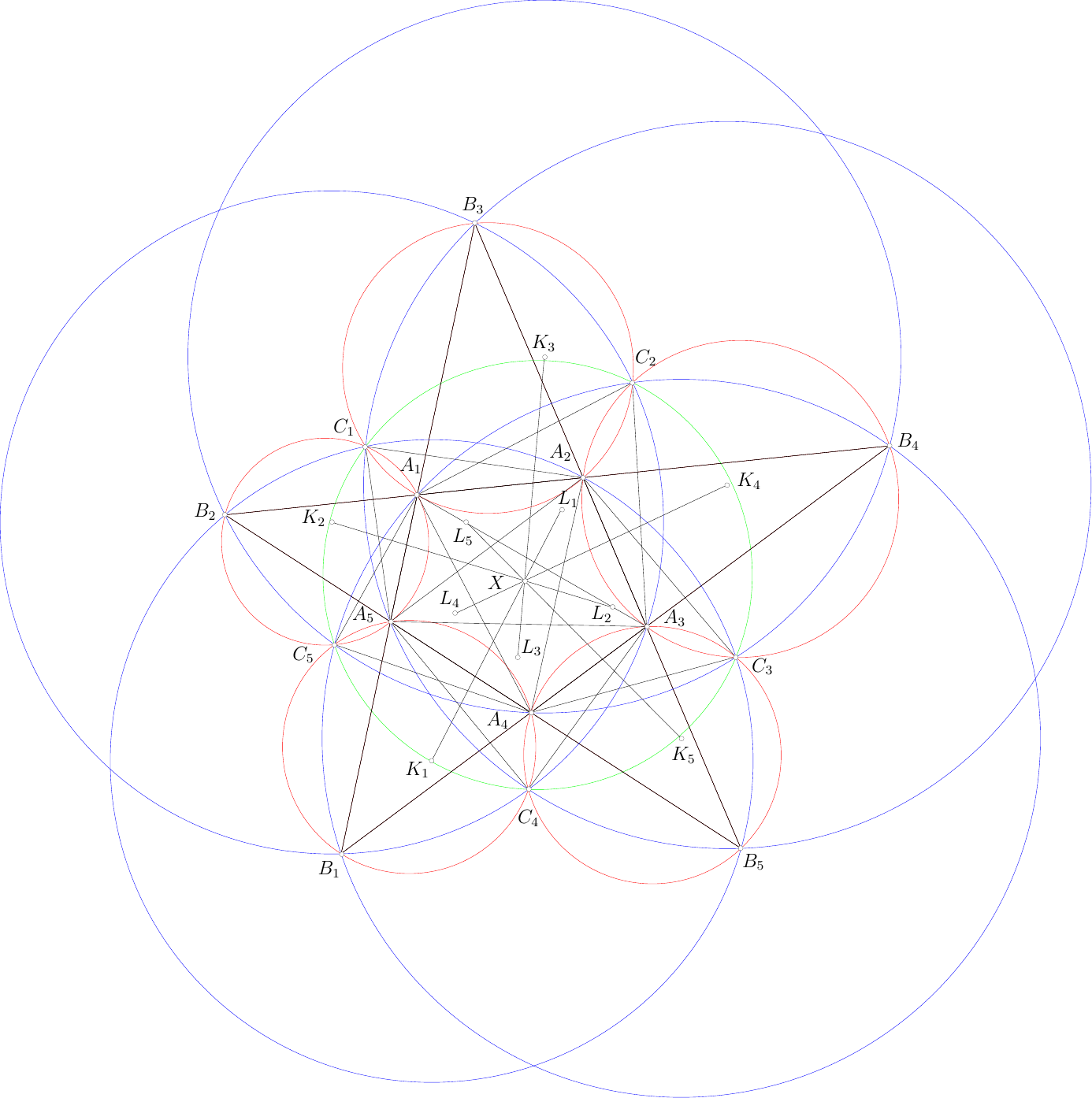}
\caption{The dual circumcenter pencil. The labels $K_i,L_i$ in this figure denote $U_i,V_i$, respectively.}\label{fig:dual}
\end{figure}

Extending the sides of a convex pentagon forms five ear triangles. Miquel's theorem says that the second intersections of consecutive ear circles are concyclic. Here we consider their circumcenters. Join each ear center to the center of the circle through the corresponding Miquel point and two outer vertices, as in Figure~\ref{fig:pencil}. The five joins form a pencil. Exchanging the pentagon and its pentagram gives another pencil, shown in Figure~\ref{fig:dual}. Both allow a parallel case, determined entirely by the side directions.

For the proof, we express the ten circumcenters of triangles bounded by the five side lines in terms of one interpolation polynomial. The resulting function of two line parameters is symmetric and affine in each variable. This allows us to combine the centers with signed area weights and obtain the pencil point. Cyclic vertices and cyclic outer vertices require further arguments. In both cases an involution pairs points on two circles, but only the outer vertex case identifies its center directly with the pencil point. The vertex case needs an identity between circle powers.

Miquel's pentagram theorem appeared in \cite{Miquel1838}. Halbeisen, Hungerb\"uhler and Loureiro \cite{HHLPentagon2023} prove it by cross ratios in Miquelian M\"obius planes. Our formulas for circumcenters use Steiner's theorem for four lines \cite{Steiner1828} and belong to the metric study of line configurations developed by Morley \cite{Morley1900}. Under the common circle hypothesis in Theorem~\ref{thm:iteration}, we obtain a correspondence with the classical branch of Morley's five circles theorem \cite{Morley1907,Dantzig1916}. The two branches of that theorem are studied in \cite{HHLTwin2023}.

Radical axes provide another link with earlier work. The incidence theorem of Fisher, Hoehn and Schr\"oder \cite{FHS2014}, extended in \cite{FSS2015}, proves that five radical axes are concurrent or parallel. It implies the theorem attributed to Takada in \cite{Shimada}. Hungerb\"uhler \cite{HungerSiblings2023} obtains further incidences from a common radical center, and Nicollier \cite{Nicollier2016} uses inversion to prove related collinearities of centers. These methods also enter our cyclic and iterated constructions.

The five circumcenter results were stated without proofs in the 2019 version of this preprint. We give the proofs here, together with formulas for the finite pencil points and the parallel cases. The notation and preliminary lemmas precede the main results in Sections~\ref{sec:pencil} to~\ref{sec:iteration}. Section~\ref{sec:morley} proves the Morley correspondence and discusses further questions.

\section{Classical results and notation}\label{sec:classical}
\begin{thm}[Miquel's pentagram theorem \cite{Miquel1838,HHLPentagon2023}]\label{thm:miquel}
Extend the sides of a convex pentagon and form the five ear triangles. The circumcircles of adjacent ears meet at a vertex of the pentagon and at a second point. The five second points are concyclic.
\end{thm}
The form needed below also allows nonconvex pentagons and crossing sides. Provided the ear circles and second intersections exist, the second points lie on a circle or a line. This is the five line form of Miquel's theorem, proved by directed angles or by the cross ratio argument in \cite{HHLPentagon2023}. We fix its indexing in Section~\ref{sec:config}. For the four line theorem, see \cite{Zhao2010}.
\begin{figure}[htbp]
\centering\includegraphics[width=.64\linewidth]{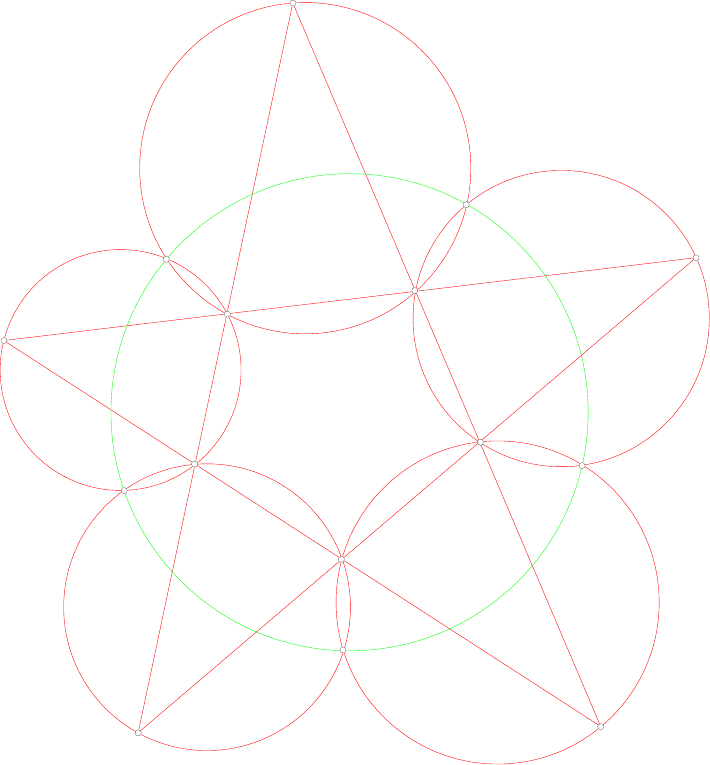}
\caption{Miquel's pentagram theorem.}\label{fig:miquel}
\end{figure}

\begin{thm}[Morley's theorem on five circles, classical branch \cite{Morley1907,Dantzig1916,HHLTwin2023}]\label{thm:five}
Let $Z_1,\ldots,Z_5$ be distinct points of a circle, normalized by a similarity to $|z|=1$, and define $P_i$ by
\[
p_i=-\frac{z_iz_{i+2}z_{i-2}}{z_{i+1}z_{i-1}}.
\]
Here the indices are taken modulo $5$. Let $\sigma_j$ be the circle with center $Z_j$ through $P_{j+2}$. Suppose that $\sigma_{i-2}\cap\sigma_{i+2}=\{P_i,Q_i\}$ with $Q_i\ne P_i$, and put $m_i=Q_{i-2}Q_{i+2}$ and $R_i=m_{i-1}\cap m_{i+1}$. Then $R_i$ lies on $\sigma_i$, and $Z_i,P_i,R_i$ are collinear.
\end{thm}
To check that $\sigma_j$ contains $P_{j-2}$ as well, put $\nu_j=z_{j+2}z_{j+4}/(z_{j+1}z_{j+3})$. Then $p_{j+2}-z_j=-z_j(1+\nu_j)$ and $p_{j-2}-z_j=-z_j(1+\nu_j^{-1})$ have the same modulus $|1+\nu_j|$, since $|\nu_j|=1$. The minus sign in the definition of $p_i$ is essential here. Reversing it gives the hidden twin of \cite{HHLTwin2023}, with different incidences.
\begin{figure}[htbp]
\centering\includegraphics[width=.64\linewidth]{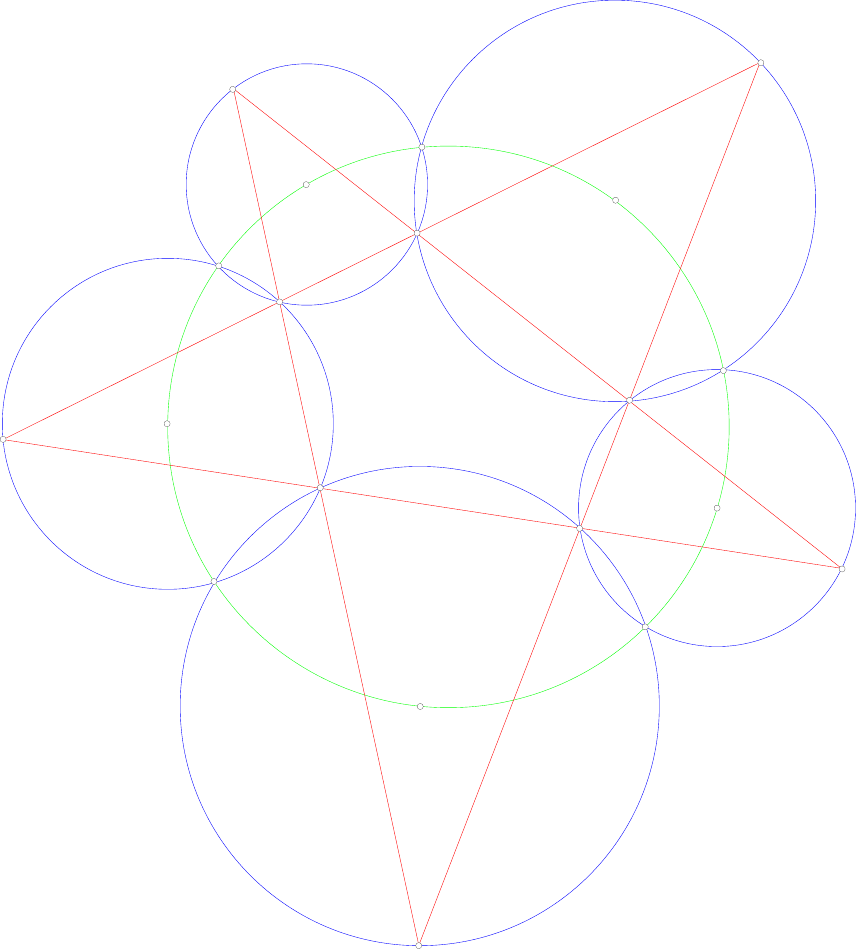}
\caption{The classical branch of Morley's theorem on five circles.}\label{fig:five}
\end{figure}

\begin{samepage}
\begin{thm}[Takada's theorem \cite{Shimada}; compare \cite{FHS2014,FSS2015}]\label{thm:takada}
Let $P_1,\ldots,P_5$ be a cyclic pentagon, with indices modulo $5$, and set
\[
R_i=P_iP_{i+2}\cap P_{i+1}P_{i-1}, \gamma_i=(P_i,P_{i+1},R_i).
\]
Let $T_i$ be the second intersection of $\gamma_{i-1}$ and $\gamma_i$, the first one being $P_i$. The points $T_i$ lie on a circle or on a line.
\end{thm}
\end{samepage}
We use the attribution in \cite{Shimada}. To prove the result, apply \cite[Theorem~1]{FSS2015}, originating in \cite{FHS2014}, with $A_i=P_i$, $B_{i,i+1}=R_i$ and $c_{i,i+1}=\gamma_i$. Its auxiliary lines are then $P_{i-1}P_{i+1}$ and $P_iP_{i+2}$, and its conclusion is that the radical axes $P_iT_i$ form a pencil.

Suppose the axes meet at $Y$, of common power $\tau$ in the five circles. For $\tau\ne0$, signed inversion about $Y$ with power $\tau$ exchanges $P_i,T_i$, since each $\gamma_i$ is invariant. The vertex circle is sent to a circle through the $T_i$, unless it passes through $Y$, in which case its image is a line. For $\tau=0$, all the $\gamma_i$ contain $Y$. But $Y\ne P_i$, as $\gamma_{i+2}$ meets the vertex circle only at $P_{i+2},P_{i+3}$. Consequently all $T_i$ coincide with $Y$.

If the axes are parallel, reflect in the line containing the five circle centers. The reflection exchanges the two intersections $P_i,T_i$ of each pair of consecutive circles and sends the vertex circle to a circle through all the $T_i$.
\begin{figure}[htbp]
\centering\includegraphics[width=.64\linewidth]{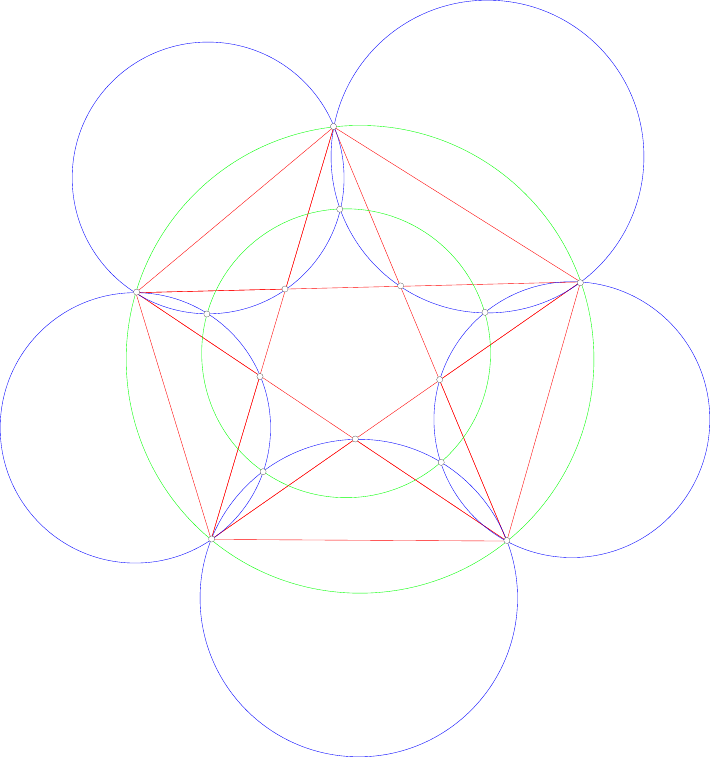}
\caption{Takada's construction for a cyclic pentagon.}\label{fig:takada}
\end{figure}

\subsection{The configuration of five lines}\label{sec:config}
From now on all indices are taken modulo $5$. Put $\ell_i=A_iA_{i+1}$ and
\begin{equation}\label{eq:configuration}
B_{i+3}=\ell_i\cap\ell_{i+2}, \omega_{i+2}=(A_i,A_{i+1},B_{i+2}),
\end{equation}
so that $\omega_{i+2}$ is the circumcircle of the ear triangle on the side $A_iA_{i+1}$. Let $K_i$ be the center of $\omega_i$, let $C_{i+1}$ be the second intersection of $\omega_{i+2}$ and $\omega_{i+3}$ (the first one is $A_{i+1}$), and let $L_i$ be the center of
\[
\lambda_i=(C_{i+1},B_{i+2},B_{i+3}).
\]
The point $C_{i+1}$ is the Miquel point of the four lines other than $\ell_{i+3}$. The five points $C_i$ lie on a circle or a line. In the former case we call their circle the Miquel circle. Miquel's theorem for four lines also shows that
\begin{equation}\label{eq:miquel-membership}
A_{i-2},A_{i-1},B_i,C_{i-2},C_{i-1}\in\omega_i
\end{equation}
and that the five points
\begin{equation}\label{eq:auxiliary}
B_{i+1},\ B_{i+4},\ A_{i+1},\ C_i,\ C_{i+2}
\end{equation}
lie on one circle, which we call $\Gamma_i$.

Any three of the five lines bound a triangle. We write $O_{jk}$ for the circumcenter of the triangle bounded by the three lines other than $\ell_j$ and $\ell_k$. The two families of centers in our theorems are
\[
K_i=O_{i,i+1},\ U_i=O_{i+2,i+4},
\]
obtained by omitting two adjacent or two nonadjacent sides. Here $U_i$ is the center of $\Gamma_i$, since $B_{i+1},B_{i+4},A_{i+1}$ are the vertices of the triangle bounded by $\ell_i,\ell_{i+1},\ell_{i+3}$. Table~\ref{tab:objects} lists the objects used later.
\begin{table}[htbp]
\centering
\renewcommand{\arraystretch}{1.15}
\begin{tabular}{c|l}
object & definition\\\hline
$\omega_i$, $K_i$ & circle $(A_{i-2},A_{i-1},B_i)$ and its center, $K_i=O_{i,i+1}$\\
$C_{i+1}$ & Miquel point of the four lines other than $\ell_{i+3}$\\
$\lambda_i$, $L_i$ & circle $(C_{i+1},B_{i+2},B_{i+3})$ and its center\\
$\Gamma_i$, $U_i$ & circle $(B_{i+1},B_{i+4},A_{i+1},C_i,C_{i+2})$ and its center, $U_i=O_{i+2,i+4}$\\
$V_i$ & center of the circle $(A_i,C_{i+1},A_{i+2})$
\end{tabular}
\caption{Objects attached to the configuration of five lines.}\label{tab:objects}
\end{table}

A configuration will be called admissible if the side directions are distinct, no three sides are concurrent, and all the specified circles, second intersections and joins are defined. In particular, a second intersection must differ from the first, a join must have distinct endpoints, and points used to define a circle center must determine a genuine circle uniquely. The notation $X$ for a pencil point always refers to a unique finite common point. All subsequent configurations are assumed admissible. Limits are included whenever they satisfy the same requirements.

\section{Complex and inversive tools}\label{sec:tools}
Lowercase letters denote complex coordinates of the corresponding points, for example $k_i$ for $K_i$. For a line we use the equation
\begin{equation}\label{eq:line}
z+t\bar z=h, |t|=1,
\end{equation}
with direction parameter $t=-e^{2\mathrm i\theta}$, where $\theta$ is a direction angle. Two lines have the same parameter precisely when they are parallel. We use the signed area convention $[a,b,c]=\frac12\operatorname{Im}\bigl(\overline{(b-a)}(c-a)\bigr)$.

\begin{lem}[Circumcenter for three lines]\label{lem:center}
Consider three lines \eqref{eq:line} with distinct parameters $r,s,t$, and let $uT^2+vT+w$ be the quadratic polynomial that takes the values $h_r,h_s,h_t$ at $r,s,t$. The triangle bounded by the three lines has circumcenter $w$ and circumradius $|u|$, and $\bar w=v+u(r+s+t)$. Moreover, if $z_{rs}$ denotes the intersection of the lines with parameters $r$ and $s$, then
\[
[z_{rs},z_{st},z_{tr}]=|u|^2\,[r^{-1},s^{-1},t^{-1}].
\]
\end{lem}
\begin{proof}
At the intersection of the lines with parameters $r,s$, subtracting their equations yields $\bar z_{rs}=u(r+s)+v$. Substitution into either equation gives $z_{rs}=w-urs$. The analogous formulas for the other two intersections put all three vertices at distance $|u|$ from $w$.

The conjugate of $z_{rs}=w-urs$, compared with $\bar z_{rs}=u(r+s)+v$ for two pairs of parameters, gives $\bar u=rst\,u$ and $\bar w=v+u(r+s+t)$. Also $rs=rst\cdot t^{-1}$, with analogous expressions for $st,tr$. Thus $z\mapsto w-urst\,z$ is the direct similarity taking $t^{-1},r^{-1},s^{-1}$ to the three vertices. Areas scale by $|u|^2$, and the cyclic permutation of the parameters preserves the sign.
\end{proof}
Thus the ratio of signed area to squared circumradius depends only on the three line directions. Its absolute value is twice the product of the sines of the triangle's angles.

\begin{lem}[Steiner's center circle \cite{Steiner1828}]\label{lem:fourcenters}
The circumcenters of the four triangles formed by four lines in general position, together with the Miquel point of the four lines, lie on one circle.
\end{lem}
\begin{proof}
Interpolate the four line constants by $AT^3+BT^2+CT+D$, taking value $h_j$ at $t_j$. Write $s_1=\sum t_j$, $s_2=\sum_{j<k}t_jt_k$, $s_4=\prod t_j$ and $G=B+As_1$. Omitting the $j$th line leaves the quadratic remainder modulo $\prod_{k\ne j}(T-t_k)$. Its coefficients are
\[
u=G-At_j, v=C-A(s_2-s_1t_j+t_j^2), w=D+As_4/t_j .
\]
By Lemma~\ref{lem:center} the triangle has center $k_j=D+As_4/t_j$ and radius $r_j=|G-At_j|$, and
\[
\bar k_j=v+u(s_1-t_j)=F-Gt_j, F=C+Bs_1+A(s_1^2-s_2).
\]
Conjugating the formula for $k_j$ also gives $\bar k_j=\bar D+\overline{As_4}\,t_j$. Since these affine expressions agree at four distinct values of $t_j$, their coefficients agree, giving $F=\bar D$ and $G=-\overline{As_4}$. Hence $|G|=|A|$. Neither is zero, since otherwise all four radii $r_j$ would vanish. Set $M=D+A^2s_4/G$. We have
\[
M-k_j=\frac{As_4(At_j-G)}{Gt_j},
\]
so $|M-k_j|=|G-At_j|=r_j$, placing $M$ on all four circumcircles. The remaining distances, $|k_j-D|=|A|$ and $|M-D|=|A|^2/|G|=|A|$, put $M$ and the four centers on the circle of radius $|A|$ about $D$. Steiner states this as assertion~2 in \cite{Steiner1828}; see also \cite[Theorem~3.1(ii)]{HungerPohle2020}.
\end{proof}
The Miquel point of four lines in general position differs from every intersection of two of the lines. Otherwise one component circumcircle would contain three distinct points of a single line. The proof also shows that the four centers $k_j=D+As_4/t_j$ are distinct, since $A\ne0$ and the $t_j$ are distinct, and that $M$ differs from each of them, since $r_j>0$.

\begin{lem}[Odd cycles of concyclic pairs]\label{lem:odd}
Let $\Omega$ and $\Omega'$ be distinct circles and let $n\ge3$ be odd. Let $P_i\in\Omega$ be distinct points and let $Q_i\in\Omega'$, with indices modulo $n$. Suppose that $P_i,Q_i,P_{i+1},Q_{i+1}$ are four distinct points on a circle or a line for every $i$. Then a signed inversion or a line reflection interchanges every pair $P_i,Q_i$.
\end{lem}
\begin{proof}
First exclude marked points from $\Omega\cap\Omega'$. If $P_i$ were an intersection point, the two circles through the adjacent quadruples would both be $\Omega'$. This would make $P_{i-1}$ and $P_{i+1}$ equal to the other intersection point, contrary to their distinctness. If $Q_i$ were an intersection point, the adjacent circles would instead be $\Omega$, forcing the $Q_i$ to alternate between the two intersection points. Oddness of $n$ rules out that alternation. For tangent circles, the same arguments force two consecutive marked points to coincide.

With the intersection points excluded, a M\"obius transformation gives one of three normal forms. The marked points remain finite in each case.

If $\Omega$ and $\Omega'$ are disjoint, we may take them concentric with radii $\rho\ne\rho'$ and write $P_i=\rho u_i$, $Q_i=\rho'v_i$ with $|u_i|=|v_i|=1$. The two chords $P_iP_{i+1}$ and $Q_iQ_{i+1}$ have the same perpendicular bisector. For a circular quadruple, both bisectors pass through the origin and the quadruple's circle center, which cannot be the origin because $\rho\ne\rho'$. For a collinear quadruple, they both pass through the origin perpendicular to its line. The bisector of the chord $\rho u,\rho u'$ has direction $\sqrt{uu'}$, so $u_iu_{i+1}=v_iv_{i+1}$. Thus the ratios $\tau_i=v_i/u_i$ satisfy $\tau_i\tau_{i+1}=1$. Going around the odd cycle forces them all to be a constant $\tau\in\{1,-1\}$, and $z\mapsto\tau\rho\rho'/\bar z$ exchanges the pairs.

If $\Omega$ and $\Omega'$ meet in two points, send these points to $0$ and $\infty$. The circles become two lines through $0$. With signed coordinates $p_i,q_i$ on these lines, the power of $0$ gives $p_ip_{i+1}=q_iq_{i+1}$. As before $q_i/p_i$ is a constant $\pm1$, and the reflection in one of the two angle bisectors swaps the pairs.

If $\Omega$ and $\Omega'$ are tangent, send the point of contact to $\infty$. Then $\Omega,\Omega'$ become parallel lines, concyclicity gives $p_i+p_{i+1}=q_i+q_{i+1}$, and since $n$ is odd we get $p_i=q_i$ for all $i$. The reflection in the midline swaps the pairs.

Returning to the original plane, we obtain an antiholomorphic M\"obius involution $\Phi$. If it fixes infinity, it is a reflection. Otherwise $\Phi(z)=(a\bar z+b)/(c\bar z+d)$, with $c\ne0$, and the involution condition at infinity gives $a/c=-\overline{d/c}$. Set $p=a/c$ to rewrite it as
\begin{equation}\label{eq:inversion}
\Phi(z)=p+\frac{\kappa}{\bar z-\bar p}.
\end{equation}
Substituting into $\Phi\circ\Phi=\mathrm{id}$ shows that $\kappa$ is real and nonzero. This is inversion with center $p$ and signed power $\kappa$, which may be negative.
\end{proof}
If a circle with center $k$ and radius $r$ does not pass through $p$, its image under \eqref{eq:inversion} has center
\begin{equation}\label{eq:image-center}
k'=p+\frac{\kappa(k-p)}{|k-p|^2-r^2}.
\end{equation}
We will need two consequences of this formula. First, the original and image centers lie on a line through $p$. Second, an invariant circle on which $\Phi$ exchanges distinct points has power $\kappa$ at $p$. For $k\ne p$, setting $k'=k$ gives $|k-p|^2-r^2=\kappa$, which is also sufficient for invariance. For $k=p$, the radius is preserved exactly when $|\kappa|=r^2$. The positive choice fixes the circle pointwise, so exchanging distinct points requires $\kappa=-r^2$, again its power at $p$.

We will also use this involution at limits where an intermediate quadruple ceases to have four distinct points. Three distinct source points with three distinct images determine an antiholomorphic M\"obius map uniquely. Its coefficients then have a continuous limit with nonzero determinant, and the identity $\Phi\circ\Phi=\mathrm{id}$ still holds. These three pairs suffice to retain the involution at such a limit.

\section{The principal circumcenter pencil}\label{sec:pencil}
For distinct $j,k$ let $a<b<c$ be the other three indices in $\{1,\ldots,5\}$. Let $T_{jk}$ be the triangle with vertices $\ell_a\cap\ell_b$, $\ell_b\cap\ell_c$, $\ell_c\cap\ell_a$, so that $O_{jk}$ is its circumcenter, and let $R_{jk}$ be its circumradius. Put
\begin{equation}\label{eq:weights}
w_{jk}=\frac{[T_{jk}]}{R_{jk}^2}.
\end{equation}
As noted after Lemma~\ref{lem:center}, $w_{jk}$ depends only on the directions of the side lines, and $|w_{jk}|$ is twice the product of the sines of the angles of $T_{jk}$.

\begin{thm}[Principal circumcenter pencil]\label{thm:pencil}
In the configuration \eqref{eq:configuration}, assume that the five lines $K_iL_i$ are defined, and put
\[
\begin{aligned}
\Sigma&=2\sum_{i=1}^5 w_{i,i+1}+\sum_{i=1}^5 w_{i+2,i+4},\\
\Xi&=2\sum_{i=1}^5 w_{i,i+1}K_i+\sum_{i=1}^5 w_{i+2,i+4}U_i .
\end{aligned}
\]
In the expression for $\Xi$, points denote their position vectors with respect to a fixed origin. If $\Sigma\ne0$, the five lines $K_iL_i$ pass through the point $X=\Xi/\Sigma$. If $\Sigma=0$, then $\Xi\ne0$ and the five lines are parallel to the vector $\Xi$.
\end{thm}
The common point is an affine combination of the ten circumcenters $O_{jk}$, with weights determined by the directions of the sides.
\begin{proof}
Let $\ell_j$ have parameter $t_j$ and constant $h_j$, and let
\[
\mathcal P(T)=\alpha T^4+\beta T^3+\gamma T^2+\delta T+\varepsilon
\]
take the values $h_j$ at $t_j$. Put $s_1=\sum t_j$, $T_0=\prod t_j$ and $h_*=\beta+\alpha s_1$. If $r,s,t$ are the parameters of the three lines other than $\ell_j,\ell_k$, the remainder of $\mathcal P$ modulo $(T-r)(T-s)(T-t)$ has constant term $\varepsilon+\bigl(\beta+\alpha(r+s+t)\bigr)rst$. By Lemma~\ref{lem:center},
\begin{equation}\label{eq:ten-centers}
o_{jk}=\varepsilon-\alpha T_0\left(t_j^{-1}+t_k^{-1}\right)+\frac{h_*T_0}{t_jt_k}.
\end{equation}
Let $\zeta_j=1/t_j$ and
\[
\mathcal O(a,b)=\varepsilon-\alpha T_0(a+b)+h_*T_0\,ab .
\]
Then $o_{jk}=\mathcal O(\zeta_j,\zeta_k)$. The points $\zeta_j$ are distinct points of the unit circle, and $\mathcal O$ is symmetric and affine in each variable. Lemma~\ref{lem:center} also gives
\begin{equation}\label{eq:weight-zeta}
w_{jk}=[\zeta_a,\zeta_b,\zeta_c].
\end{equation}

Here $a<b<c$ are the indices other than $j,k$. To express $L_i$ in the same coordinates, consider the four lines $\ell_{i-1},\ell_i,\ell_{i+1},\ell_{i+2}$ whose Miquel point is $C_{i+1}$. The circumcircles centered at $O_{i+3,i}$ and $O_{i+3,i+2}$ share the chord $C_{i+1}B_{i+2}$. The pair centered at $O_{i+3,i+1}$ and $O_{i+3,i+4}$ shares $C_{i+1}B_{i+3}$. Their lines of centers are the perpendicular bisectors of these two chords, so their intersection is $L_i$.

The map $z\mapsto\mathcal O(\zeta_{i+3},z)$ sends $\zeta_k$ to $O_{i+3,k}$. It is a direct similarity, since its linear coefficient is nonzero. A zero coefficient would make the four circumcircles concentric; as they share $C_{i+1}$, they would then coincide, giving three intersections with $\ell_i$. Applying the similarity to the two lines just described gives
\begin{equation}\label{eq:L-affine}
\begin{aligned}
l_i&=\mathcal O(\zeta_{i+3},Y_i),\\
Y_i&=\zeta_i\zeta_{i+2}\cap\zeta_{i+1}\zeta_{i+4}.
\end{aligned}
\end{equation}

Besides $[a,b,c]+[a,c,d]=[a,b,d]+[b,c,d]$ we use the determinant identity
\begin{equation}\label{eq:area2}
[a,b,c]d+[a,c,d]b=[a,b,d]c+[b,c,d]a .
\end{equation}
Triangulate the labelled pentagon $\zeta_1\cdots\zeta_5$ by diagonals. Let $S$ be the sum, over the three triangles, of the signed area of the triangle times the value of $\mathcal O$ at the two remaining vertices, and let $W$ be the sum of the three signed areas. Changing a diagonal in $abcd$ changes its contribution to $S$ from $[a,b,c]\mathcal O(d,e)+[a,c,d]\mathcal O(b,e)$ to $[a,b,d]\mathcal O(c,e)+[b,c,d]\mathcal O(a,e)$, where $e$ is the fifth vertex. The two area identities make these equal because $\mathcal O(\cdot,e)$ is affine. They also preserve $W$. Any triangulation can be reached by diagonal flips, so the sums are independent of the triangulation. This argument uses only determinant identities and does not require the labels to follow the order around the circle.

For a fixed $i$, abbreviate $(\zeta_i,\ldots,\zeta_{i+4})$ by $(a,b,c,d,e)$. Choosing the three triangles with common vertex $c$, we have
\begin{equation}\label{eq:fan}
\begin{aligned}
S&=[c,d,e]\,\mathcal O(a,b)+[c,e,a]\,\mathcal O(b,d)+[c,a,b]\,\mathcal O(d,e),\\
W&=[c,d,e]+[c,e,a]+[c,a,b].
\end{aligned}
\end{equation}
Here $\mathcal O(a,b)=k_i$, and since $Y_i=ac\cap be$,
\[
[c,e,a]\,b+[c,a,b]\,e=\bigl([c,e,a]+[c,a,b]\bigr)Y_i .
\]
Applying the affine map $\mathcal O(d,\cdot)$ and using \eqref{eq:L-affine}, we obtain
\begin{equation}\label{eq:area-pencil}
S=[c,d,e]\,k_i+\bigl(W-[c,d,e]\bigr)\,l_i .
\end{equation}
The points $c,d,e$ are distinct points of the unit circle, so $[c,d,e]\ne0$. If $W\ne0$, then $S/W$ is a real affine combination of $k_i$ and $l_i$ for every $i$, so all five lines $K_iL_i$ pass through $S/W$. If $W=0$, then $S=[c,d,e](k_i-l_i)$, which is nonzero because $K_iL_i$ is defined, and every line $K_iL_i$ is parallel to $S$.

We can now recover the stated weights. The three values of $\mathcal O$ in \eqref{eq:fan} are $K_i=O_{i,i+1}$, $U_{i-1}=O_{i+1,i+3}$ and $K_{i+3}=O_{i+3,i+4}$. Formula~\eqref{eq:weight-zeta} identifies their area coefficients as $w_{i,i+1}$, $w_{i+1,i+3}$ and $w_{i+3,i+4}$; cyclic shifts preserve the signs. Hence
\begin{equation}\label{eq:fan-centers}
\begin{aligned}
S&=w_{i,i+1}K_i+w_{i+1,i+3}U_{i-1}+w_{i+3,i+4}K_{i+3},\\
W&=w_{i,i+1}+w_{i+1,i+3}+w_{i+3,i+4}.
\end{aligned}
\end{equation}
Sum these expressions over the five choices of $i$. An adjacent pair $\{j,j+1\}$ is counted twice and a nonadjacent pair once, giving $5S=\Xi$ and $5W=\Sigma$.
\end{proof}
\begin{rem}[Directions of the sides]
Let $\theta_j$ be a direction angle of $\ell_j$. Then $\zeta_j=-e^{-2\mathrm i\theta_j}$, and the shoelace formula gives
\[
\Sigma=5W=\tfrac52\operatorname{Im}\sum_j\bar\zeta_j\zeta_{j+1}=-\tfrac52\sum_j\sin2(\theta_{j+1}-\theta_j).
\]
For five distinct joining lines, parallelism is therefore equivalent to $\sum_j\sin2(\theta_{j+1}-\theta_j)=0$, or to $\sum_j\sin2A_j=0$ for a convex pentagon with angles $A_1,\ldots,A_5$. Moving the sides while retaining their directions leaves this condition unchanged. If $\Sigma\ne0$, normalize the three weights in \eqref{eq:fan-centers} by their sum. They express $X$ relative to $K_i,U_{i-1},K_{i+3}$ and give its barycentric coordinates when those centers are noncollinear.
\end{rem}

\section{The dual pencil}\label{sec:dual}
\begin{samepage}
\begin{thm}[Dual circumcenter pencil]\label{thm:dual}
Let $U_i$ be the center of $\Gamma_i$ and $V_i$ the center of $(A_i,C_{i+1},A_{i+2})$. If the five lines $U_iV_i$ are defined, they pass through one point or are parallel.
\end{thm}
\begin{proof}
Apply Theorem~\ref{thm:pencil} to the pentagon $\widetilde A_j=B_{2j-1}$. Its sides and outer vertices are
\[
\widetilde\ell_j=\ell_{2j-2}, \widetilde B_j=A_{2j-3}.
\]
With this indexing, the ear circle $(B_{2j},B_{2j+2},A_{2j+2})$ is $\Gamma_{2j+1}$ and has center $U_{2j+1}$. Consecutive ear circles meet at $\widetilde A_j=B_{2j-1}$ and at $\widetilde C_j=C_{2j}$. The circle corresponding to $\lambda_j$ is consequently defined by
\[
\widetilde C_{j+1},\ \widetilde B_{j+2},\ \widetilde B_{j+3}
\ =\ C_{2j+2},\ A_{2j+1},\ A_{2j+3},
\]
so its center is $V_{2j+1}$. Theorem~\ref{thm:pencil} now applies to the five joins $U_iV_i$, with $i=3,5,2,4,1$.
\end{proof}
\end{samepage}
\FloatBarrier
\begin{rem}
Thus the dual pencil is the principal pencil for the ordering $\ell_5,\ell_2,\ell_4,\ell_1,\ell_3$. Its point is given by the same weighted formula, with $K_i,U_i$ interchanged and the complementary triples in \eqref{eq:weights} ordered accordingly. The direction calculation gives $\sum_j\sin2(\theta_{j+2}-\theta_j)=0$ as the parallel condition when the five joins $U_iV_i$ are distinct.

There are $4!/2=12$ cyclic orders of five lines, counting an order and its reverse as the same. Each gives a circumcenter construction. The relabelling in the proof pairs them into six pentagon and pentagram constructions, although different orders can give the same pencil in special cases.
\end{rem}

\section{Cyclicity and the line of centers}\label{sec:cyclic}
The next two proofs are first given for configurations with distinct intermediate points, and in Theorem~\ref{thm:cyclic} for the inversion case. Here is the justification for passing to limits. Normalize cyclic vertices as $a_j=e^{\mathrm i\theta_j}$; for cyclic outer vertices, use the same parameterization on the dual pentagon of Theorem~\ref{thm:dual}. The constructed coordinates are rational in $a_j,\bar a_j$. Clearing denominators turns the excluded coincidences and the reflection condition into polynomial equations. None holds identically, as the regular pentagon on the unit circle shows. In that configuration,
\[
c_i=\varphi a_i, b_i=-\varphi^2a_{i+1}, \varphi=\frac{1+\sqrt5}{2}.
\]
Here the circles of the $A_i$, $B_i$ and $C_i$ are distinct, all the needed quadruples have four distinct points, and the two involutions are signed inversions. In the cyclic vertex case, \eqref{eq:cyclic-system} also has nonzero determinant. Thus none of the excluded polynomial conditions holds identically. Their zero sets have empty interior, leaving a dense set of configurations on which to prove the results.

At an admissible limit, the defining triples remain noncollinear, so their circle centers vary continuously. The same holds for $X=\Xi/\Sigma$, since uniqueness and finiteness imply $\Sigma\ne0$ by Theorem~\ref{thm:pencil}. Collinearity passes to the limit. The involution does too, by the observation after Lemma~\ref{lem:odd}, since its image circle supplies three distinct image points.

\begin{thm}[Cyclic outer vertices]\label{thm:outer}
Suppose that the $B_i$ lie on a circle with center $O$, and let $J$ be the center of the Miquel circle. If the pencil $K_iL_i$ has a unique finite point $X$, then $O,J,X$ are collinear. Moreover, $X$ is the center of the signed inversion that swaps $B_i$ and $C_{i+1}$ for every $i$.
\end{thm}
\begin{proof}
The circle of the $B_i$ differs from the Miquel circle. If they were equal, $\Gamma_i$ would meet that circle at $B_{i+1},B_{i+4}$, since $A_{i+1}$ does not lie on it. Indeed, otherwise the line $\ell_i$ would meet the circle of the $B_i$ at the three points $A_{i+1}$, $B_{i+1}$ and $B_{i+3}$. Applying this to $\Gamma_i$ and $\Gamma_{i-2}$ would give $C_i=B_{i+4}$ for every $i$, hence $C_{i+1}=B_i$. This contradicts the observation following Lemma~\ref{lem:fourcenters}, because $B_i$ is a vertex of the four lines defining $C_{i+1}$.

Pair $P_i=B_i$ with $Q_i=C_{i+1}$. By \eqref{eq:auxiliary}, the circle $\Gamma_i$ contains the pairs with indices $i+1$ and $i+4$, namely $P_{i+1}=B_{i+1}$, $Q_{i+1}=C_{i+2}$, $P_{i+4}=B_{i+4}$ and $Q_{i+4}=C_i$. Reading the pairs in the order $1,4,2,5,3$ gives the odd cycle required by Lemma~\ref{lem:odd}. Its involution $\Phi$ exchanges $B_i$ and $C_{i+1}$, and persists at admissible limits by the argument above.

The three points $B_i,C_{i-2},C_{i-1}$ on $\omega_i$ are sent to $C_{i+1},B_{i+2},B_{i+3}$, the defining points of $\lambda_i$. Thus $\Phi(\omega_i)=\lambda_i$. A reflection would put every $K_iL_i$ perpendicular to its axis, contradicting the existence of a unique finite pencil point. Hence $\Phi$ is a signed inversion. Denote its center by $P$. The centers $K_i,L_i$ and $P$ are collinear by \eqref{eq:image-center}. Hence $P$ is the common pencil point $X$. The inversion also sends the circle of the $B_i$ to the Miquel circle, putting their centers $O,J$ on a line through this same point.
\end{proof}
Nicollier \cite{Nicollier2016} uses inversion for related collinearities of centers. The pairing needed here is $B_i\leftrightarrow C_{i+1}$.
\begin{figure}[htbp]
\centering\includegraphics[width=.98\linewidth]{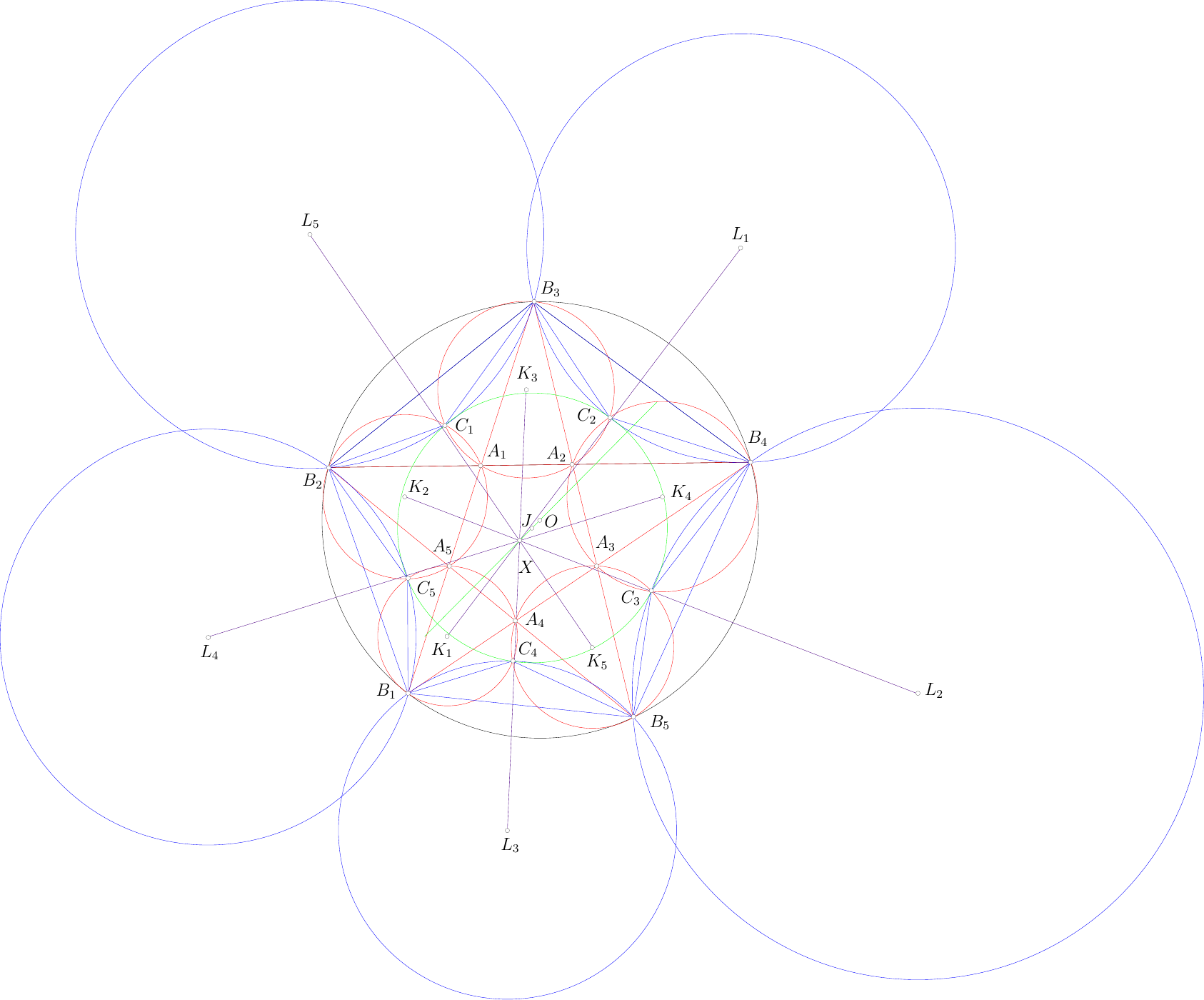}
\caption{Theorem~\ref{thm:outer}, with the outer vertices $B_i$ on a circle.}\label{fig:outer}
\end{figure}
\FloatBarrier

\begin{thm}[Cyclic vertices]\label{thm:cyclic}
Suppose that the $A_i$ lie on a circle $(O)$, and let $J$ be the center of the Miquel circle. If the pencil $K_iL_i$ has a unique finite point $X$, then $O,J,X$ are collinear. More precisely, if the involution of Lemma~\ref{lem:odd} that swaps $A_i$ and $C_i$ is a signed inversion with center $P$ and power $\kappa$, then
\[
\overrightarrow{OX}=\frac{\Pow\bigl(P,(O)\bigr)}{\kappa}\,\overrightarrow{OJ}.
\]
\end{thm}
\begin{proof}
Normalize $(O)$ to the unit circle and write $a_i$ for $A_i$. The vertex circle differs from the Miquel circle. Otherwise $\omega_i$, which is not the vertex circle, would force $C_{i-2}=A_{i-1}$ and $C_{i-1}=A_{i-2}$. The circle $\omega_{i+1}$ would also force $C_{i-1}=A_i$, a contradiction. Here $\omega_i$ differs from the vertex circle because its outer vertex $B_i$ does not lie on that circle.

Apply Lemma~\ref{lem:odd} to the pairs $A_i,C_i$, using the circles $\omega_i$ in \eqref{eq:miquel-membership} for consecutive pairs. The resulting involution persists at admissible limits. In the inversion case, write it as \eqref{eq:inversion} and put
\[
q=|p|^2, N=\kappa-q .
\]
For a point $a_i$ of the unit circle, \eqref{eq:inversion} gives
\begin{equation}\label{eq:C-inversion}
c_i=\frac{p+Na_i}{1-\bar p a_i}.
\end{equation}
The inversion preserves $\omega_i$ and exchanges distinct points on it, so its power at $p$ is $\kappa$ by the criterion after \eqref{eq:image-center}. In terms of its power function
\begin{equation}\label{eq:circle-power}
F_i(z)=|z|^2-\bar k_i z-k_i\bar z+n_i,
\end{equation}
we get $F_i(p)=\kappa$, that is, $n_i=\bar k_ip+k_i\bar p+N$.

Fix $i$ and put $(a,b,c,d,e)=(a_i,a_{i+1},a_{i+2},a_{i+3},a_{i+4})$. For $\omega_i$, the three bounding chords are $cd$, $de$, $ea$. Their intersection and center equations give
\[
\bar k_i=\frac{a-c}{ea-cd},\ k_i=de\,\bar k_i,\ n_i=\frac{da-ce}{ea-cd},\ F_i(b)=-\frac{(a-c)(b-d)(b-e)}{b(ea-cd)}.
\]
The circle $\lambda_i$ passes through
\[
u=B_{i+2}=\frac{ea(b+c)-bc(e+a)}{ea-bc},\ v=B_{i+3}=\frac{ab(c+d)-cd(a+b)}{ab-cd},\ m=\Phi(b)=C_{i+1}.
\]
To find the center of $\lambda_i$, introduce the real quadratic function
\[
\mathcal H(z)=\kappa|z|^2+(N+1)(\bar p z+p\bar z)+2N^2+N-q .
\]
The calculation below establishes
\begin{equation}\label{eq:local-identity}
\mathcal H(z)F_i(b)=(N+1)^2F_i(z)\text{ for }z=u,v,m.
\end{equation}

At $m=\Phi(b)$, formula~\eqref{eq:C-inversion}, with $d_b=|b-p|^2$, yields $d_b\mathcal H(m)=\kappa(N+1)^2$. Invariance of $\omega_i$ gives $d_bF_i(m)=\kappa F_i(b)$. Eliminating $d_b$ proves \eqref{eq:local-identity} at $m$.

At $u,v$, use the circles of the chord chains $(e,a,b,c)$, $(a,b,c,d)$ and $(c,d,e,a)$, namely $\omega_{i+2}$, $\omega_{i+3}$ and $\omega_i$. Their common power $\kappa$ at $p$ yields the linear system $\mathsf M\,(p,\bar p,N)^{\mathsf T}=y$, where
\begin{equation}\label{eq:cyclic-system}
\mathsf M=\begin{pmatrix}
c-e&ab(c-e)&bc-ea\\
d-a&bc(d-a)&cd-ab\\
a-c&de(a-c)&ea-cd
\end{pmatrix},\ y=\begin{pmatrix}ac-eb\\bd-ac\\da-ce\end{pmatrix}.
\end{equation}
Let $\Delta=\det\mathsf M$, and let $\Delta_j$ be the determinant of $\mathsf M$ with its $j$th column replaced by $y$. Since $y$ is the combination of the columns of $\mathsf M$ with coefficients $p,\bar p,N$, multilinearity gives
\[
\Delta_1=p\Delta, \Delta_2=\bar p\Delta, \Delta_3=N\Delta .
\]
Adding the two determinant expansions for $\Lambda=\Delta+\Delta_3$ and factoring gives
\[
\Lambda=-(a-c)^2(a-d)(b-d)(b-e)(c-e).
\]
All factors are nonzero, since the vertices are distinct. The equality $\Lambda=(N+1)\Delta$ then proves $\Delta\ne0$ and $N+1\ne0$. Using $\kappa=N+p\bar p$ and $q=p\bar p$, we clear the denominators in $\mathcal H(z)$ to get
\[
\Delta^2\mathcal H(z)=(\Delta_3\Delta+\Delta_1\Delta_2)\,z\bar z+\Lambda(\Delta_2z+\Delta_1\bar z)+2\Delta_3^2+\Delta_3\Delta-\Delta_1\Delta_2 ,
\]
which is polynomial in $a,\ldots,e$, $z$ and $\bar z$. At $u,v$, use the following conjugates and abbreviations.
\[
\begin{array}{c|c|c}
z&\bar z&R_z\\\hline
u&\dfrac{a-b-c+e}{ea-bc}&\dfrac{ab(c-e)^2}{(ea-bc)^2}\\[2ex]
v&\dfrac{a+b-c-d}{ab-cd}&\dfrac{bc(a-d)^2}{(ab-cd)^2}
\end{array}
\]
Substitution gives
\[
\begin{aligned}
F_i(z)&=F_i(b)\,R_z,\\
\Delta^2\mathcal H(z)&=\Lambda^2R_z .
\end{aligned}
\]
Since $\Lambda^2=(N+1)^2\Delta^2$, these two identities give \eqref{eq:local-identity} for $z=u$ and $z=v$.

Set $\mu_i=(N+1)^2/F_i(b)$, which is real and defined by the displayed factorization of $F_i(b)$. If $\mu_i=\kappa$, the function $\mathcal H-\kappa F_i$ is affine linear and has value
\[
\mathcal H(p)-\kappa F_i(p)=\kappa(N+1)\ne0
\]
at $p$, whereas \eqref{eq:local-identity} would make this function vanish at the noncollinear points $u,v,m$. This contradiction proves $\mu_i\ne\kappa$. Hence $\mathcal H-\mu_iF_i=0$ is a circle equation, and its three points $u,v,m$ identify it as $\lambda_i$. Dividing the coefficient of $\bar z$ by the negative of the coefficient of $|z|^2$ gives the center
\[
l_i=\frac{\mu_ik_i+(N+1)p}{\mu_i-\kappa}.
\]
For $x=-(N+1)p/\kappa$ this becomes
\[
l_i-x=\frac{\mu_i}{\mu_i-\kappa}\,(k_i-x),
\]
whose real coefficient puts $x$ on $K_iL_i$ for every $i$. Thus $x=X$. The image of the unit circle under $\Phi$ is the genuine Miquel circle, so the inversion center is off the unit circle and $q\ne1$. Its image center is given by \eqref{eq:image-center} as
\[
j=\frac{N+1}{1-q}\,p,\ x=\frac{q-1}{\kappa}\,j .
\]
Here $O=0$ and $q-1=\Pow(P,(O))$, giving the claimed vector relation.

It remains to include the reflection case in the collinearity conclusion. A reflection preserving the $\omega_i$ has their centers on its axis. Dividing each row of \eqref{eq:cyclic-system} by its third entry gives rows $(\bar k,k,1)$ for three of those centers, so the matrix is singular. The reflection case thus belongs to the exceptional set $\Delta=0$. The density argument at the start of this section and continuity of $O,J,X$ extend the collinearity to this case.
\end{proof}
Since $N+1\ne0$ in the inversion case, $X=O$ holds exactly when $P=O$. In general, however, the inversion center $P$ and the pencil point $X$ are different.
\begin{figure}[htbp]
\centering\includegraphics[width=.98\linewidth]{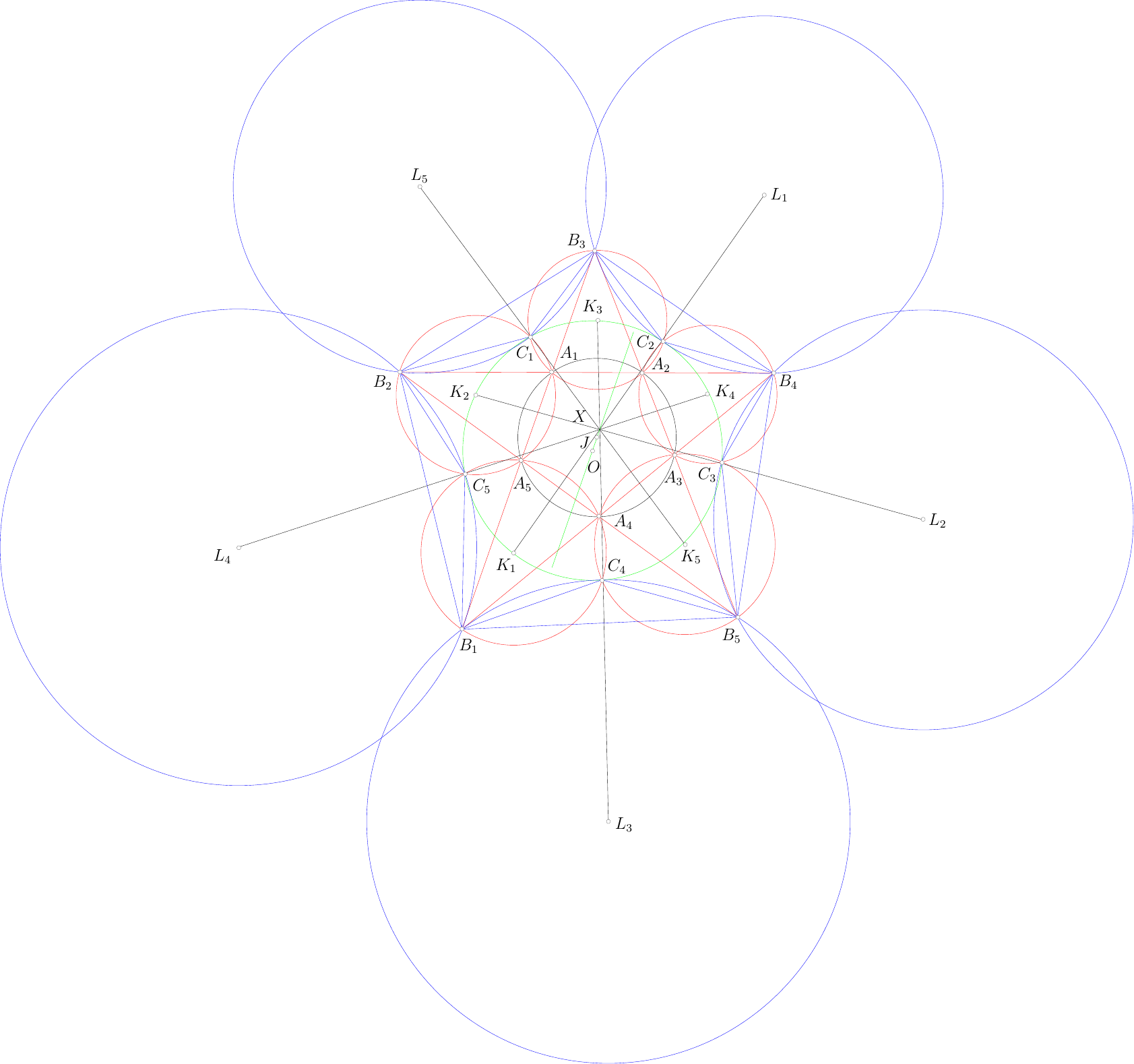}
\caption{Theorem~\ref{thm:cyclic}, with the vertices $A_i$ on a circle.}\label{fig:cyclic}
\end{figure}
\FloatBarrier

\section{An iterated Miquel construction}\label{sec:iteration}
\begin{thm}[Iteration on a common circle]\label{thm:iteration}
Assume that $K_1,\ldots,K_5$ and $C_1,\ldots,C_5$ lie on one circle $(O)$. Define
\begin{equation}\label{eq:D-def}
D_i=K_iK_{i+2}\cap K_{i+1}K_{i+3},
\end{equation}
and let $E_i$ be the second intersection of the circles $(D_i,D_{i+1},K_{i+2})$ and $(D_{i+1},D_{i+2},K_{i+3})$, the first one being $D_{i+1}$. The points $E_i$ lie on a circle or a line, and $K_i,L_i,E_i$ are collinear for every $i$. Assume in addition that the $E_i$ determine a genuine circle with center $J$ and that the pencil $K_iL_i$ has a unique finite point $X$. Then a signed inversion with center $X$ interchanges $K_i,E_i$ for every $i$, and $O,J,X$ are collinear.
\end{thm}
\begin{proof}
For each $j$, Lemma~\ref{lem:fourcenters} puts the four centers $O_{jk}$ $(k\ne j)$ and the Miquel point of the four lines other than $\ell_j$ on one circle. Two of these centers are $K_j$ and $K_{j-1}$, and the Miquel point is $C_{j-2}$. These three points are distinct by the observation following Lemma~\ref{lem:fourcenters}. All three lie on $(O)$, so this circle is $(O)$. Hence all ten centers $O_{jk}$ lie on $(O)$.

By \eqref{eq:ten-centers},
\[
o_{jk}=\bigl(\varepsilon-\alpha T_0/t_j\bigr)+T_0\bigl(h_*/t_j-\alpha\bigr)/t_k ,
\]
so the four centers with fixed $j$ lie on a circle about $\varepsilon-\alpha T_0/t_j$. They are distinct and already belong to $(O)$, making this its center for every $j$. Since the $t_j$ are distinct, we must have $\alpha=0$ and $(O)$ centered at $\varepsilon$.

Normalize $(O)$ to the unit circle. Then $o_{jk}=\beta T_0/(t_jt_k)$ with $|\beta T_0|=1$. Taking a square root $\xi^2=\beta T_0$ and writing $u_j=\xi/t_j$ gives
\begin{equation}\label{eq:center-products}
|u_j|=1, o_{jk}=u_ju_k, k_i=u_iu_{i+1}.
\end{equation}
The $u_j$ are distinct, so consecutive $K_i$ are distinct. Equality of two nonconsecutive $K_i$ would make one of the lines in \eqref{eq:D-def} undefined. Thus all five $K_i$ are distinct. The argument that gave \eqref{eq:L-affine}, with the rotation $z\mapsto u_{i+3}z$ in place of $\mathcal O(\zeta_{i+3},\cdot)$, now gives $l_i=u_{i+3}Y_i$ with $Y_i=u_iu_{i+2}\cap u_{i+1}u_{i+4}$.

By \eqref{eq:D-def}, $D_iD_{i+1}=K_{i+1}K_{i+3}$. Thus the pentagon $D_1\cdots D_5$ has the side lines $K_{i+1}K_{i+3}$, and in the notation \eqref{eq:configuration} its outer vertices are exactly the points $K_i$. Its ear circles are the circles $(D_i,D_{i+1},K_{i+2})$, and the points $E_i$ are its Miquel points. By Miquel's theorem for five lines they lie on a circle or a line. Two of these chords sharing an endpoint cannot be parallel. Any two without a common endpoint meet at one of the defined points $D_i$. Thus the side directions of $D$ are distinct. Nor can three sides be concurrent. Among any three of these chords, two share an endpoint $K_i$, and the third cannot pass through that endpoint. The observation after Lemma~\ref{lem:fourcenters} now ensures $E_i\ne K_i$, so we can use the line $K_iE_i$.

Let $\Theta_i=(K_{i+1},K_{i+4},D_{i+1})$. This is the circle $\Gamma_i$ of \eqref{eq:auxiliary} for the pentagon $D_1\cdots D_5$, so it contains $E_{i-1}$ and $E_{i+1}$, and $\Theta_{i-1}$, $\Theta_{i+1}$ meet at $K_i$ and $E_i$. Let $v_i$ be the center of $\Theta_i$. The triangle $K_{i+1}K_{i+4}D_{i+1}$ is bounded by three chords of the unit circle, and \eqref{eq:center-products} gives
\[
v_i=u_iu_{i+1}\,\frac{u_{i+2}-u_{i+4}}{u_i-u_{i+1}}.
\]
For $(a,b,c,d,e)=(u_i,\ldots,u_{i+4})$ this means
\[
\begin{gathered}
k_i=ab,\ l_i=d\,\frac{ac(b+e)-be(a+c)}{ac-be},\\
v_{i-1}=ea\,\frac{b-d}{e-a},\ v_{i+1}=bc\,\frac{d-a}{b-c},
\end{gathered}
\]
and over a common denominator
\begin{equation}\label{eq:perpendicular-identity}
v_{i+1}-v_{i-1}=\frac{ac-be}{(a-e)(b-c)}\,(l_i-k_i).
\end{equation}
Since $\bar a=1/a,\ldots,\bar e=1/e$, conjugation changes the sign of the factor in front of $l_i-k_i$, so this nonzero factor is purely imaginary. Hence $K_iL_i$ is perpendicular to the line of centers of $\Theta_{i-1}$ and $\Theta_{i+1}$. Their radical axis $K_iE_i$ is perpendicular to the same line and also passes through $K_i$. Therefore $K_iE_i=K_iL_i$, and $K_i,L_i,E_i$ are collinear.

Suppose now that the pencil has a unique finite point $X$. Because each $K_iL_i$ is the radical axis of $\Theta_{i-1},\Theta_{i+1}$, the point $X$ has equal powers in these two circles. Stepping by $2$ cycles through all five indices, giving a common power $\tau$ for every $\Theta_i$. Along their shared chords, the directed lengths satisfy $XK_i\cdot XE_i=\tau$.

If $\tau$ were zero, then $X=K_i$ or $X=E_i$ for each $i$. At most one $K_i$ equals $X$, so at least four $E_i$ would coincide there. The remaining points could not determine the assumed circle $(J)$. We therefore have $\tau\ne0$, and signed inversion about $X$ with power $\tau$ exchanges $K_i,E_i$. It sends $(O)$ to $(J)$, so their centers are collinear with $X$ by \eqref{eq:image-center}.
\end{proof}
\begin{figure}[htbp]
\centering\includegraphics[width=.98\linewidth]{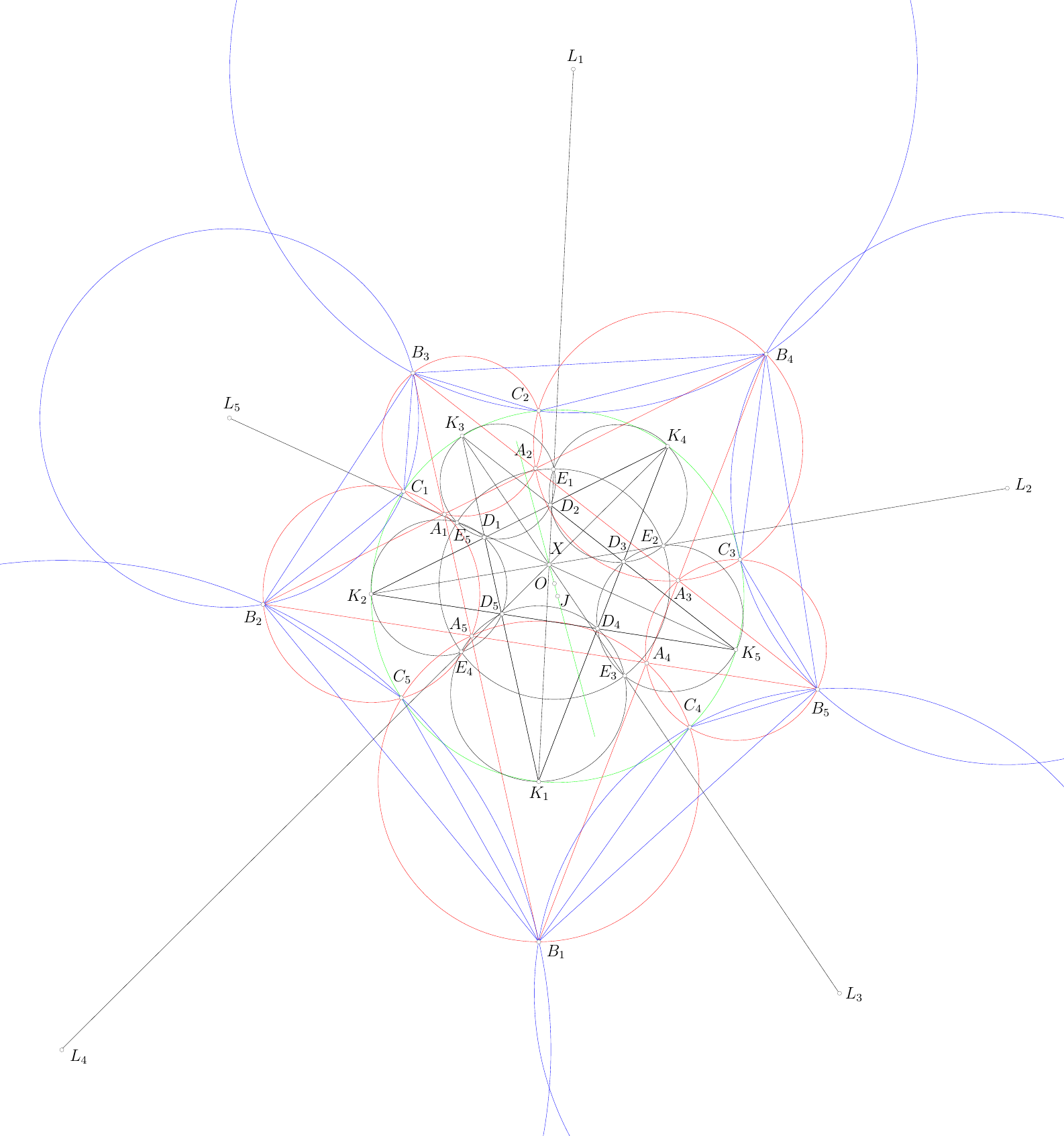}
\caption{The iterated configuration of Theorem~\ref{thm:iteration}.}\label{fig:iteration}
\end{figure}

\begin{rem}[The common circle condition]
For admissible configurations the common circle hypothesis of Theorem~\ref{thm:iteration} is equivalent to $\alpha=0$ in the polynomial $\mathcal P$ of the proof of Theorem~\ref{thm:pencil}. One direction was shown above. Conversely, if $\alpha=0$, then \eqref{eq:ten-centers} gives $|o_{jk}-\varepsilon|=|\beta|$. Moreover, for the four lines other than $\ell_j$ the cubic of Lemma~\ref{lem:fourcenters} is $\mathcal P$ itself, so the Miquel point of these lines lies on the circle with center $\varepsilon$ and radius $|\beta|$ as well.
\end{rem}
\begin{samepage}
\section{Morley's theorem and further questions}\label{sec:morley}
\subsection{The common circle case as a Morley configuration}
In the situation of Theorem~\ref{thm:iteration} the configuration is an instance of Theorem~\ref{thm:five}. Normalize $(O)$ to the unit circle and use the following notation.
\[
\begin{array}{c|ccccc}
\text{Theorem~\ref{thm:five}} & Z_i & \sigma_i & P_i & Q_i & R_i\\ \hline
\text{this paper} & K_i & \omega_i & C_{i+1} & A_{i+1} & B_i
\end{array}
\]
To check the correspondence, use \eqref{eq:miquel-membership} to obtain $\omega_{i-2}\cap\omega_{i+2}=\{C_{i+1},A_{i+1}\}$. The line $m_i=Q_{i-2}Q_{i+2}$ is then $A_{i-1}A_{i-2}=\ell_{i-2}$, and its prescribed intersection is $m_{i-1}\cap m_{i+1}=\ell_{i-3}\cap\ell_{i-1}=B_i$.

\end{samepage}
\begin{samepage}
\begin{prop}\label{prop:morley-link}
Under the hypothesis of Theorem~\ref{thm:iteration}, the correspondence above gives the classical branch of Theorem~\ref{thm:five}. Consequently
\[
c_{i+1}=-\frac{k_ik_{i+2}k_{i-2}}{k_{i+1}k_{i-1}}.
\]
The points $B_i,K_i,C_{i+1}$ are collinear. Moreover, the center $U_i$ of $\Gamma_i$ lies on $(O)$, and its coordinate is
\[
\frac{k_{i-1}k_{i+1}}{k_i}.
\]
\end{prop}
\end{samepage}
\begin{proof}
The proof of Theorem~\ref{thm:iteration} gives $\alpha=\varepsilon=0$, so $\mathcal P(T)=\beta T^3+\gamma T^2+\delta T$, with $|\beta|=1$ because $|\beta T_0|=1$, and $k_i=u_iu_{i+1}$ in the notation \eqref{eq:center-products}. Put $s_1=\sum_j t_j$. Each line is real, so $\mathcal P(t_j)=t_j\overline{\mathcal P(t_j)}$. Multiplying by $t_j^2$ and using $\bar t_j=1/t_j$, we see that every $t_j$ is a root of
\[
\beta T^5+\gamma T^4+\delta T^3-\bar\delta T^2-\bar\gamma T-\bar\beta ,
\]
which is therefore $\beta\prod_j(T-t_j)$. Comparing the coefficients of $T^4$ gives $\gamma=-\beta s_1$.

The point $C_{i+1}$ is the Miquel point of the four lines other than $\ell_m$, where $m=i+3$. For these four lines the cubic in Lemma~\ref{lem:fourcenters} is $\mathcal P$ itself, so $A=\beta$, $D=0$, $s_4=T_0/t_m$ and $G=\gamma+\beta(s_1-t_m)=-\beta t_m$. Lemma~\ref{lem:fourcenters} gives
\[
c_{i+1}=\frac{\beta^2T_0/t_m}{-\beta t_m}=-\frac{\beta T_0}{t_m^2}=-u_m^2
=-\frac{k_ik_{i+2}k_{i-2}}{k_{i+1}k_{i-1}} .
\]
This is the sign of the classical branch, and Theorem~\ref{thm:five} shows that $B_i,K_i,C_{i+1}$ are collinear.

Finally, $\Gamma_i$ is the circumcircle of the triangle bounded by $\ell_i,\ell_{i+1},\ell_{i+3}$, so by \eqref{eq:center-products} its center has coordinate
\[
o_{i+2,i+4}=u_{i+2}u_{i+4}=\frac{k_{i-1}k_{i+1}}{k_i},
\]
which has modulus $1$.
This is also the center of the extra circle in \cite[Theorem~5]{HHLTwin2023}.
\end{proof}
\subsection{Further questions}
Can the two circumcenter pencils be proved synthetically? The proof of Theorem~\ref{thm:pencil} uses signed areas only after the centers have been expressed by complex interpolation. A geometric derivation of those area relations would remove the interpolation step. Another route would be to construct circles whose consecutive radical axes are $K_iL_i$ or $U_iV_i$. Theorem~\ref{thm:iteration} does this under the common circle hypothesis, while the constructions in \cite{HungerSiblings2023} suggest looking for such circles more generally.

Two steps in the cyclic constructions still depend on algebra. In Theorem~\ref{thm:cyclic}, the determinant calculation establishes \eqref{eq:local-identity}, relating the power function of $\omega_i$ to an equation of $\lambda_i$. Can a coaxal family explain this relation? In Theorem~\ref{thm:iteration}, it is the perpendicularity in \eqref{eq:perpendicular-identity} that needs a geometric proof. The radical axis argument would then finish the collinearity, and inversion would give the line of centers. The Morley correspondence supplies incidences that could be used in seeking such a proof.

There are also twelve pencils attached to the cyclic orders of the same five lines. They share five Miquel points, each determined by an omitted line. Which relations among their finite pencil points hold for arbitrary admissible configurations?

With more than five lines, one would need closure conditions for the new constructions. The circle incidence and closing results in \cite{FSS2015,HungerClosing2025} offer examples to compare. For Miquelian planes, there is an earlier obstacle. Circumcenters are Euclidean objects and are not preserved by M\"obius transformations, although the Miquel incidences are preserved. An extension would first need to replace the centers by objects defined in that setting.

For related circumcenter questions in cyclic hexagons, see \cite{Dergiades2014}. The discussions \cite{DaoDiscussion,CohlDiscussion,TranHexagon,TranCyclic,TranBicentric} record further configurations that motivated these questions.

\section*{Acknowledgements}
The author thanks Professor Norbert Hungerb\"uhler for his valuable comments and suggestions, which helped to improve and complete this paper. The author also thanks Alexander Skutin for reading the manuscript and Hiroshi Okumura for information about Takada's theorem.

\end{document}